\documentclass[12pt,fleqn]{article}

\usepackage[utf8]{inputenc}
\usepackage[T1]{fontenc}

\usepackage{url}
\usepackage{makecell}
\usepackage{amsmath}
\usepackage{amsthm}

\usepackage{graphicx}
\usepackage{epstopdf}
\usepackage{authblk}

\usepackage[a4paper,margin=1in]{geometry}

\newcommand{\eqb}{\begin{equation}}
\newcommand{\eqe}{\end{equation}}
\newcommand{\al}{\alpha}

\title{Designing sustainable barn-type houses: Optimal shapes for minimal envelope and energy use}

\author[1]{Ewa Rokita-Magdziarz}
\author[2]{Barbara Gronostajska}
\author[3]{Marcin Magdziarz}
\affil[1]{Rokita-Project Architectural Office, Forsycjowa 7, Wrocław, 51-253, Poland}
\affil[2]{Department of Architectural and Construction Design, Faculty of Architecture, Wrocław University of Science and Technology, Wyspianskiego 27, Wrocław, 50-370, Poland}
\affil[3]{Hugo Steinhaus Center, Department of Applied Mathematics, Faculty of Pure and Applied Mathematics, Wrocław University of Science and Technology, Wyspianskiego 27, Wroclaw, Poland}
\affil[*]{Correspondence should be addressed to Marcin Magdziarz at \texttt{marcin.magdziarz@pwr.edu.pl}}

\date{} 

\begin{document}

\maketitle


\begin{abstract}
Barn-type houses have become one of the most popular single-family housing
typologies in Poland and across Europe due to their simplicity, functionality,
and potential for energy efficiency.
Despite their widespread use, systematic
methods for optimizing their geometry in terms of envelope area and energy performance remain limited.
This paper develops a rigorous mathematical
framework for determining the optimal proportions of barn-type houses with
respect to minimizing the external surface area while satisfying constraints of
either fixed volume or fixed floor area. Closed-form solutions for the optimal
width, length, and height are derived as explicit functions of the roof slope,
together with formulas for the minimal achievable surface. A recently introduced
dimensionless compactness measure is also calculated, allowing quantitative
assessment of how far a given design deviates from the theoretical optimum.
The methodology is applied to case studies of three existing houses, showing
that while some designs deviate significantly from optimal compactness,
others already closely approximate it. The results confirm that theoretical
optimization can lead to meaningful reductions in construction costs and energy demand. To support practical
implementation, two original freely available software tools were developed, enabling architects and engineers to perform
optimization analyses.
\end{abstract}

\thispagestyle{empty}

\section{Introduction}

The built environment is widely recognized as one of the most significant contributors to global energy demand, greenhouse gas (GHG) emissions, and natural resource consumption. Buildings account for over half of worldwide energy use and more than one-third of solid waste flows, while also constituting the single greatest source of anthropogenic GHG emissions \cite{UN2014,IEA2018}. This problem is exacerbated by rapid urbanization: between now and 2030, the global urban population will increase by an additional one billion people, which translates into the construction of the equivalent of a new city of 1.5 million inhabitants every week for more than a decade \cite{Angel2012}. Without fundamental changes in design strategies, such growth is expected to double carbon emissions by mid-century \cite{IEA2018}, crossing thresholds associated with catastrophic and irreversible climate change. Thus, the efficient design of building forms is both an environmental imperative and a pressing scientific challenge \cite{Pomponi2018,LCA2017,LondonStrategy2018,California2017}.

Among the multiple determinants of a buildings environmental performance, the geometry of the building envelope plays a critical role. The envelope directly mediates thermal exchanges between indoor and outdoor environments, thereby shaping heating and cooling demands and determining operational carbon emissions \cite{Ramesh2010}. It also represents one of the most resource-intensive components of construction, often relying on materials with high embodied energy and carbon footprints. Minimizing the external surface area for a given internal volume or floor area is therefore an intuitive and effective pathway to enhance sustainability \cite{Steemers2003,Ratti2005a,Ratti2005b}. This fundamental observation has motivated both theoretical and applied research on compactness in architecture \cite{d2019compactness}.

Barn-type houses provide an especially relevant case study in this context. In Poland, they now constitute more than 10\% of new single-family dwellings \cite{pol1}, reflecting their popularity across Europe due to their simplicity, elegance, cost efficiency, and adaptability to both modern and traditional aesthetics. Their compact volumes and clean structural lines lend themselves naturally to optimization with respect to envelope minimization, while their widespread use makes them a meaningful target for broader sustainability improvements. Yet, despite this potential, systematic mathematical approaches for determining their optimal proportions remain scarce.

As stated in \cite{d2019compactness}, research on building forms can be broadly categorized into two complementary traditions: \emph{classification-seeking} and \emph{form-seeking} approaches \cite{Steadman2000,Steadman2014,Julia2017,Hargreaves2017}. Classification-seeking work includes surveys and taxonomies of building stock, which are crucial for policy analysis and urban morphology studies. Examples include Steadman's pioneering surveys of buildings in England and Wales \cite{Steadman2000,Steadman2014}, Bayesian-based archetype modeling of residential stock in the United States \cite{Julia2017}, and density-to-form translation methods applied to the UK context \cite{Hargreaves2017}. Such efforts provide a systematic understanding of existing building typologies, thereby establishing the empirical foundation for large-scale energy policy and retrofitting strategies.

Form-seeking approaches, by contrast, aim to identify geometries that optimize performance relative to given criteria. This tradition traces back to \cite{Martin1972,March1976}. Subsequent contributions include analyses of energy use across urban built forms \cite{Steemers2003}, work on energy consumption and urban texture \cite{Ratti2005a,Ratti2005b}, and empirical demonstrations of the correlation between exposed surface area and energy use across the London building stock \cite{Steadman2000}. More recent studies highlight the influence of compactness and passive design strategies in warmer climates \cite{Vartholomaios2017}.

A wide range of optimization efforts fall under this category. For office buildings, Catalina et al.\ \cite{Catalina2011} and Ourghi et al.\ \cite{Ourghi2007} quantified the effect of form on annual energy use, while Schlueter and Thesseling \cite{Schlueter2009} developed  tools to integrate energy considerations. Multi-criteria optimization approaches considering construction costs, heating demand, and emissions were pioneered by Jedrzejuk and Marks \cite{Jedrzejuk2000,Jedrzejuk2002}. Other research has sought to integrate solar considerations: Hachem et al.\ \cite{Hachem2012,Hachem2016} investigated settlement-level housing designs for solar exploitation, Okeil \cite{Okeil2010} proposed the Residential Solar Block, and Caruso et al.\ \cite{Caruso2013} developed geometric forms minimizing direct solar gains. Jin and Jeong \cite{Jin2017} employed genetic algorithms for optimizing free-form shapes with respect to thermal loads. Non-rectangular forms have been explored for minimizing cooling loads \cite{Steadman2000}, while computational fluid dynamics models have identified building forms promoting air quality in microenvironments \cite{Cheshmehzangi2016}. Economic dimensions have also been incorporated: Chau et al.\ \cite{Chau2007} and Helal et al.\ \cite{Helal2016} analyzed the optimality of building heights under cost and revenue constraints, revisiting earlier structural economics established by Khan \cite{Khan1972}.

The consensus emerging from these diverse strands is that compact building forms those with minimized external surfaces for a given functional space systematically reduce both operational and embodied energy impacts. Building upon this insight, recent advances have introduced scale-independent compactness measures that rigorously quantify deviation from geometric optimality \cite{d2019compactness}. These theoretical developments are highly relevant for barn-type houses, whose simple yet variable proportions make them ideal for analytical optimization. By deriving closed-form expressions for optimal proportions under constraints of fixed volume and fixed floor area, and by applying compactness metrics, the present paper integrates barn-type houses into this global discourse on sustainable building forms.

In what follows, we provide a systematic treatment of the problem of
optimizing the geometry of barn-type houses with respect to minimizing the
external surface area of the building envelope.

In Section \ref{analiza} we introduce the mathematical
framework for analyzing the optimal shape of a barn-type house under the
constraint of fixed volume. We derive explicit formulas for the optimal
parameters---width, length, and height---as functions of the roof slope angle,
and we present closed-form expressions for the minimal achievable surface area.

In Section \ref{measure} we recall and apply a recently
proposed dimensionless compactness measure \cite{d2019compactness}. This measure allows us to
quantitatively evaluate how far a given design deviates from the theoretical
optimum, in a scale-independent way. We present analytical formulas and
graphical representations of the compactness measure for barn-type houses,
illustrating its use in design practice.

In Section \ref{area} we turn to a second optimization
problem that is especially relevant for architectural practice: the case of
fixed floor area. Since investors often specify the required usable floor area
together with room height and roof slope, we formulate the corresponding
optimization problem and solve it explicitly.

Section \ref{cases} presents three detailed case studies of
recently built barn-type houses. We compare their actual parameters with
the theoretically optimal ones obtained from the formulas in the previous
sections.

We also create and describe two software tools developed to
facilitate the application of our methods. The programs, freely available
online, implement the optimization procedures for both fixed volume and fixed
floor area scenarios, enabling architects and engineers to quickly assess and
improve their designs.

Finally, the last section concludes the paper.

\section{Optimal shape of the barn-type house -- analysis}
\label{analiza}

Barn-type houses are currently one of the most popular single-family houses chosen by investors. It is estimated that they account for over 10\% of all houses currently being built in Poland \cite{pol1}.
Barn-type houses combine energy and cost savings, practicality, elegance, and simplicity in a single design. Their clear, uncomplicated structure keeps construction costs low while allowing for generous, open interiors that can be tailored to both modern and traditional aesthetics. With their compact shape, they achieve excellent energy performance when combined with quality insulation and sustainable materials.

\begin{figure}[t]
\centering
\includegraphics[width=11cm]{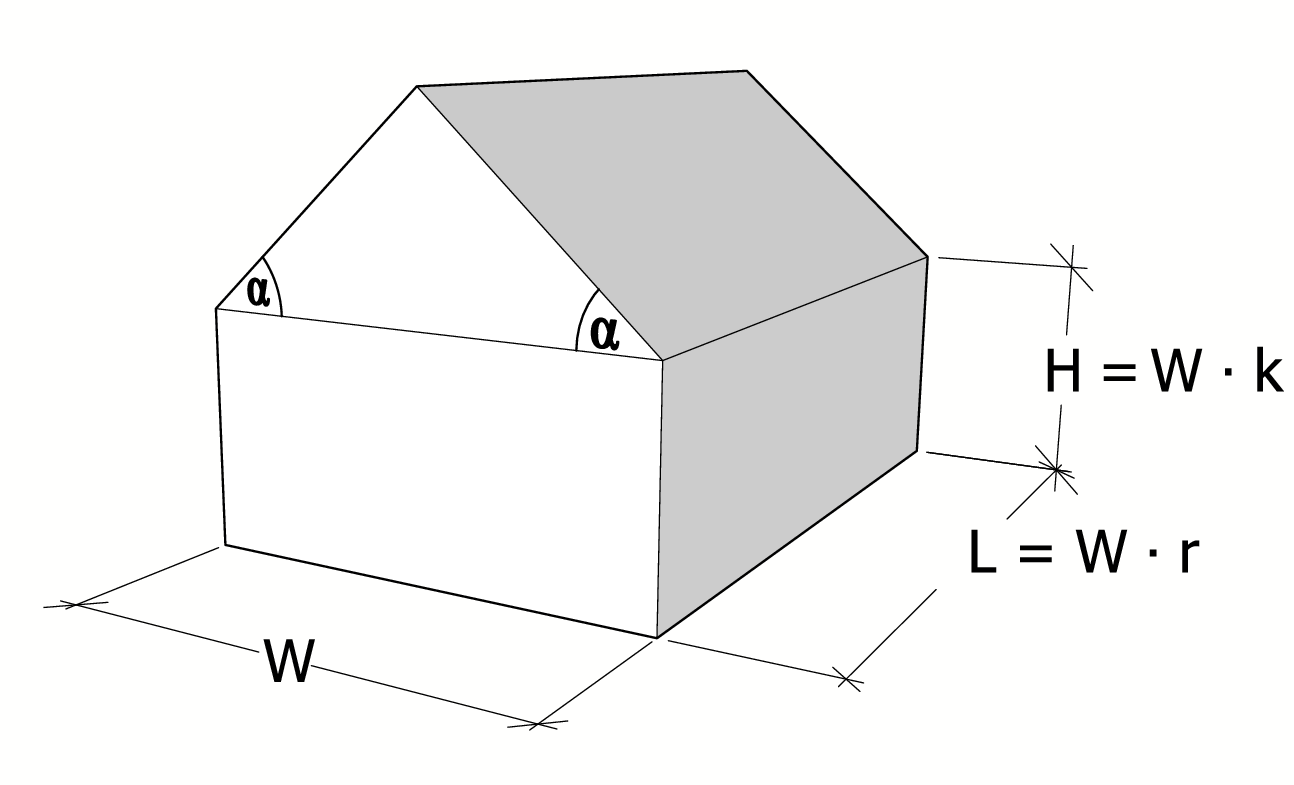}
\caption{Typical shape of the barn-type house.  }
\label{barn}
\end{figure}

A typical shape of the barn-type house is presented in Fig. \ref{barn}. It is parametrized here by:
\begin{itemize}
    \item $W$ -- width,
    \item $L$ -- length,
    \item $H$ -- height,
    \item $\alpha$ - roof slope angle.
\end{itemize}
Additionally we introduce the footprint aspect ratio $r=L/W$ and the slenderness aspect ratio $k=H/W$.

The roof slope angle $\alpha$ is assumed to be constant. It is usually specified in the local development plan.
The volume of the habitable area of the house is equal to
\eqb
\label{objetosc}
V=W L H.
\eqe
We assume for simplicity that the attic is non-habitable.

The external surface (envelope--walls and roof) of the barn-type house is given by:
\eqb
\label{powierzchnia}
S=2WH+2LH+\frac{LW}{\cos(\alpha)}+ W^2 \tan(\alpha)/2.
\eqe
Note that the part of the building in contact with the ground is not included in $S$.

Now, the main goal of this section is to find the optimal shape of the barn house. More precisely, for a given fixed volume $V$ we will find the parameters $W$, $L$ and $H$ such that the external surface of the building $S$ will be minimal. For this purpose, let us rewrite surface $S$ in \eqref{powierzchnia} as a function of two variables $r$ and $k$. Recall that $L=Wr$, $H=Wk$ and $V=WLH$, which implies that $W=(\frac{V}{rk})^{1/3}$. Moreover, $V$ and $\alpha$ are assumed fixed. Therefore, we get that
\begin{align}
S=S(r,k) &=2W^2 k + 2 W^2 rk + W^2 \frac{r}{\cos(\al)}+W^2\tan(\al)/2
\nonumber \\
&= W^2 \left(2k + 2rk + \frac{r}{\cos(\al)}+\tan(\al)/2\right) \nonumber \\
&= V^{2/3}\left( \frac{2k + 2rk + \frac{r}{\cos(\al)}+\tan(\al)/2}{(rk)^{2/3}}\right).
\label{powierzchnia2}
\end{align}
To find the minimum of the function $S(r,k)$ it is enough to find the minimum of the function in brackets in \eqref{powierzchnia2}. Let us denote it by
\eqb
\label{gamma_f}
\gamma(r,k)=\frac{2k + 2rk + \frac{r}{\cos(\al)}+\tan(\al)/2}{(rk)^{2/3}}.
\eqe
Let us calculate the partial derivatives $\frac{\partial \gamma(r,k)}{\partial r}$ and $\frac{\partial \gamma(r,k)}{\partial k}$. We have
\[
\frac{\partial \gamma(r,k)}{\partial r}=
\frac{k(r - \sin(\al) - 4k\cos(\al) + 2rk\cos(\al))}{3\cos(\al)(rk)^{5/3}},
\]
\[
\frac{\partial \gamma(r,k)}{\partial k}=
-\frac{r \left( 2r + \sin(\alpha) - 2k \cos(\alpha) - 2rk \cos(\alpha) \right)}
       {3 \cos(\alpha)  (rk)^{5/3}}.
\]
To find the minimum of the function $\gamma(r,k)$ we solve the system of equations
\begin{equation}
\label{sysytem}
\left\{
\begin{aligned}
&\frac{\partial \gamma(r,k)}{\partial r} = 0, \\[6pt]
&\frac{\partial \gamma(r,k)}{\partial k} = 0.
\end{aligned}
\right.
\end{equation}
After some tedious calculations we get that for $\al\in(0,\pi/2)$ the positive solution of \eqref{sysytem} that provides minimum of the function $\gamma(r,k)$ is given by
\eqb
\label{r_min}
r_{min} =\sqrt{\sin(\al) + 1/4} + 1/2
\eqe
\eqb
\label{k_min}
k_{min} =\frac{\sqrt{4\sin(\al) + 1} + 1}{4\cos(\al)}
\eqe
The corresponding minimal external surface of the barn-type house is equal to
\eqb
\label{S_min}
S_{min}=\frac{V^{2/3} \left( 3\sin(\alpha) + 6\sqrt{\sin(\alpha) + \tfrac{1}{4}} + 3 \right)}
     {2\cos(\alpha) \left( \dfrac{2\sin(\alpha) + 2\sqrt{\sin(\alpha) + \tfrac{1}{4}} + 1}{4\cos(\alpha)} \right)^{2/3}}.
\eqe
The remaining optimal parameters of the burn house are as follows:
\eqb
\label{W_min}
W_{min}=V^{1/3}\left( \frac{4 \cos(\alpha)}{2\sin(\alpha) + \sqrt{4\sin(\alpha) + 1} + 1} \right)^{\tfrac{1}{3}},
\eqe
\eqb
\label{L_min}
L_{min}=
V^{1/3} \left( \sqrt{\sin(\alpha) + 1/4} + 1/2 \right)
\left( \frac{4 \cos(\alpha)}{2\sin(\alpha) + \sqrt{4\sin(\alpha) + 1} + 1} \right)^{\tfrac{1}{3}},
\eqe
\eqb
\label{H_min}
H_{min}=V^{1/3}\frac{\left( \sqrt{4\sin(\alpha) + 1} + 1 \right)}{4\cos(\alpha)}
\left( \frac{4 \cos(\alpha)}{2\sin(\alpha) + \sqrt{4\sin(\alpha) + 1} + 1} \right)^{\tfrac{1}{3}}.
\eqe
Formulas \eqref{S_min}--\eqref{H_min} can be used to design a barn house with minimal external surface $S_{min}$ for a given volume $V$. This will reduce construction costs and provide reasonable energy performance of the building.

As an example, let us take $V=300\; m^3$ and $\al=\pi/6=30^\circ$. Then, from equations \eqref{r_min}--\eqref{H_min} we get the following optimal parameters:
$r_{min}=1,3660$, $k_{min}=0,7887$, $S_{min}=238,7161\;m^2 $, $W_{min}=6,5301\;m$, $L_{min}=8,9203\;m$, $H_{min}=5,1501\;m$.

It should be underlined that from a practical point of view, the parameters obtained are perfectly suited for real-life applications in house construction.

\begin{figure}[t]
\centering
\includegraphics[width=14cm]{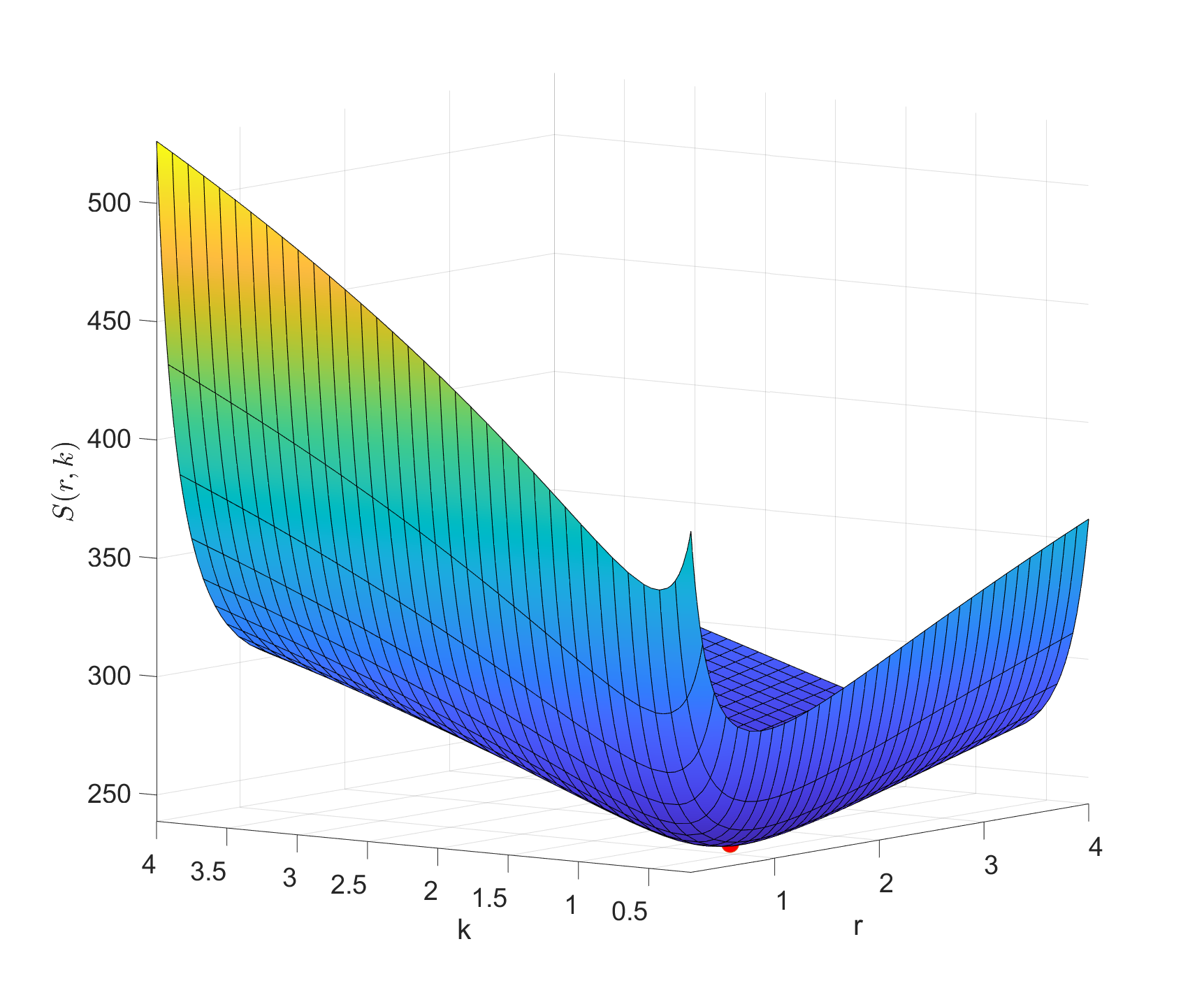}
\caption{Plot of the external surface $S(r,k)$ as a function of  footprint aspect ratio $r$ and slenderness aspect ratio $k$. Here $V=300 \; m^3$ and $\al=\pi/6=30^\circ$. Red dot is the point corresponding to the minimal surface. }
\label{S_example}
\end{figure}

In Fig. \ref{S_example} we see the graph of the external surface $S(r,k)$ as a function of footprint aspect ratio  $r$ and slenderness aspect ratio $k$ with parameters $V=300 \; m^3$ and $\al=\pi/6=30^\circ$. The minimal surface point is denoted by the red dot.

\begin{figure}[t]
\centering
\includegraphics[width=12cm]{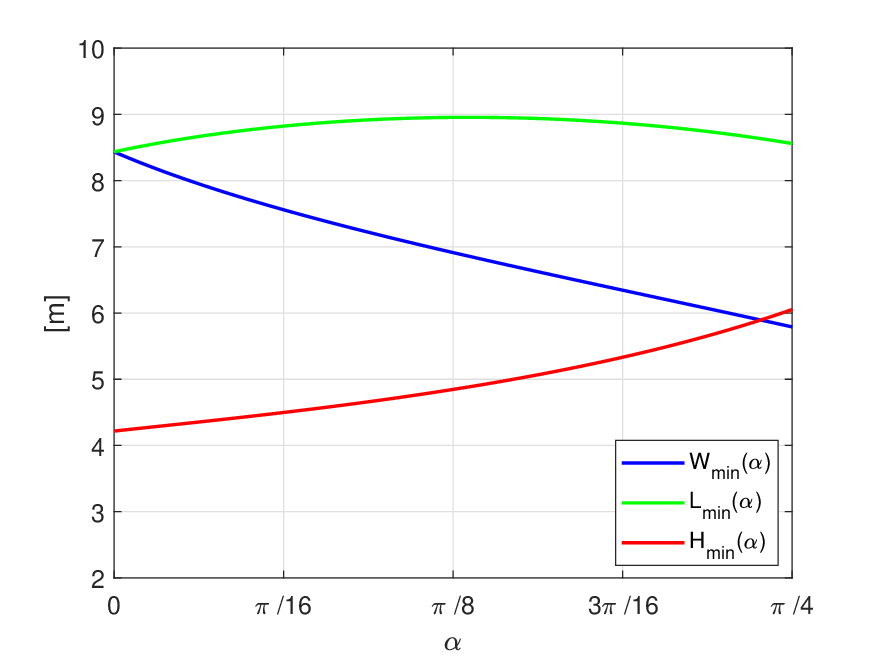}
\caption{Plots of the optimal parameters $W_{min}$, $L_{min}$ and $H_{min}$ as functions of $\al$. Here $V=300 \; m^3$.  }
\label{alpha_dependence}
\end{figure}

It is interesting to analyze how the optimal parameters $W_{min}$, $L_{min}$ and $H_{min}$ change with changing $V$ and $\al$.
From \eqref{W_min}--\eqref{H_min} we can see that in each equation there is a prefactor $V^{1/3}$, which describes the dependence of the optimal parameters on the volume $V$. We observe a more interesting situation analyzing the dependence on $\al$. In Fig. \ref{alpha_dependence} we can see plots of the optimal parameters $W_{min}$, $L_{min}$ and $H_{min}$ as functions of $\al$. Interestingly enough, $W_{min}$ decreases with increasing $\al$, $H_{min}$ increases with increasing $\al$, whereas $L_{min}$ first increases and then decreases with increasing $\al$.

\section{Compactness measure}
\label{measure}

A common way to quantify the relationship between a building's external surface $S$ and its volume $V$ is through the surface-to-volume ratio $S/V$ \cite{depecker2001design,hegger2012energy}. Lower $S/V$ values are typically preferred as they indicate a higher compactness. However, a key limitation of this ratio is that it is not scale-independent: two geometrically identical forms (e.g., cubes) of different sizes will yield different $S/V$ values. Therefore, recently a novel dimensionless compactness metric has been proposed, designed to rigorously quantify the compactness of building geometries exclusively as a function of form, thereby eliminating dependency on absolute size or volume \cite{d2019compactness}.
This compactness measure of a given building form is defined as the ratio
\eqb
\label{measure_def}
\frac{S}{S_{min}},
\eqe
where $S$ is its external surface
area of the building and $S_{min}$ denotes the minimum external surface area necessary to enclose volume $V$.
This newly developed measure can help in the early design phase by quantitatively assessing the deviation of a building form from optimal compactness and delineating the range of alternative geometries that may approach this optimum.

Let us now calculate the compactness measure for the barn-type house
with external surface $S$ and volume $V$. From \eqref{powierzchnia2} and \eqref{S_min} we find that
\begin{align}
\frac{S}{S_{min}}&=\frac{2k + 2rk + \frac{r}{\cos(\al)}+\tan(\al)/2}{(rk)^{2/3}}
\cdot \frac{ \left( \dfrac{2\sin(\alpha) + 2\sqrt{\sin(\alpha) + \tfrac{1}{4}} + 1}{4\cos(\alpha)} \right)^{2/3}\cdot 2\cos(\alpha)}{  3\sin(\alpha) + 6\sqrt{\sin(\alpha) + \tfrac{1}{4}} + 3} \nonumber \\
&= \gamma(r,k)\cdot \frac{ \left( \dfrac{2\sin(\alpha) + 2\sqrt{\sin(\alpha) + \tfrac{1}{4}} + 1}{4\cos(\alpha)} \right)^{2/3}\cdot 2\cos(\alpha)}{  3\sin(\alpha) + 6\sqrt{\sin(\alpha) + \tfrac{1}{4}} + 3}
\label{measure_formula}
\end{align}
Note that $\frac{S}{S_{min}}$ calculated above is indeed dimensionless and does not depend on $V$. It factorizes into $\gamma(r,k)$ defined in \eqref{gamma_f} and the $\al$-dependent reminder. Moreover, $\frac{S}{S_{min}}\geq 1$ with equality achieved for the optimal shape of the house.
The measure provides information on the optimality of
the shape of the barn-type house and can be used at the design
stage to verify its efficiency and to show potential change directions leading to energy and cost savings. We will use this compactness measure in Section 'Case studies' to check the optimality of certain barn-type houses that exist in reality.

\begin{figure}[t]
\centering
\includegraphics[width=14cm]{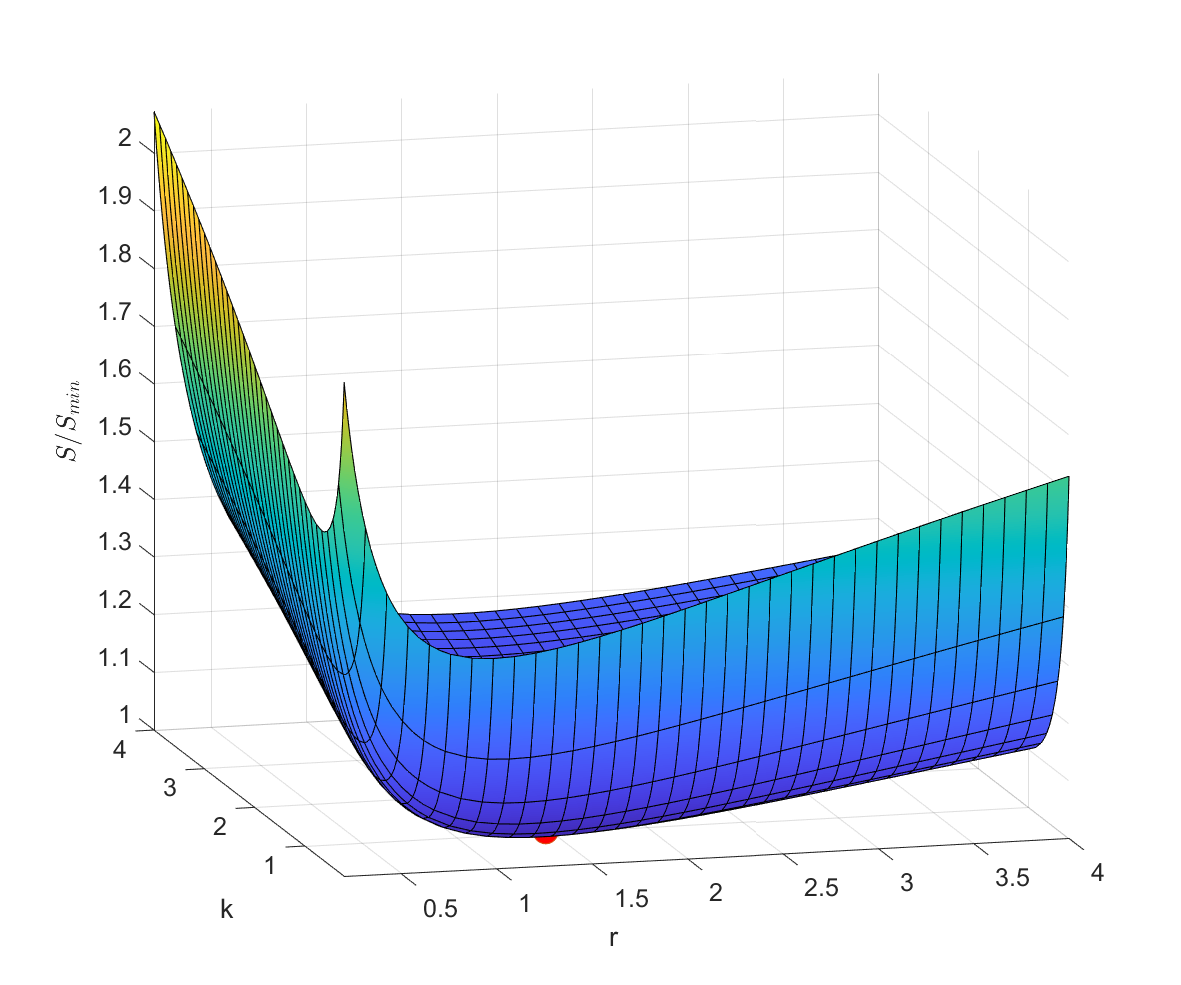}
\caption{Plot of the compactness measure $\frac{S}{S_{min}}$ as a function of  footprint aspect ratio $r$ and slenderness aspect ratio $k$. Here $\al=\pi/4=45^\circ$. Red dot is the point corresponding to the global minimum equal to 1. }
\label{measure_plot}
\end{figure}

A plot of the compactness measure $\frac{S}{S_{min}}$ as a function of $r$ and $k$ for $\al=\pi/4=45^\circ$ is shown in Fig. \ref{measure_plot}. It reaches its global miniumum equal to one for $r=r_{min}=1,4783$ and $k=k_{min}=1,0453$ (recall that $r_{min}$ and $k_{min}$ are given by \eqref{r_min} and \eqref{k_min}, respectively).

\begin{figure}[t]
\centering
\includegraphics[width=14cm]{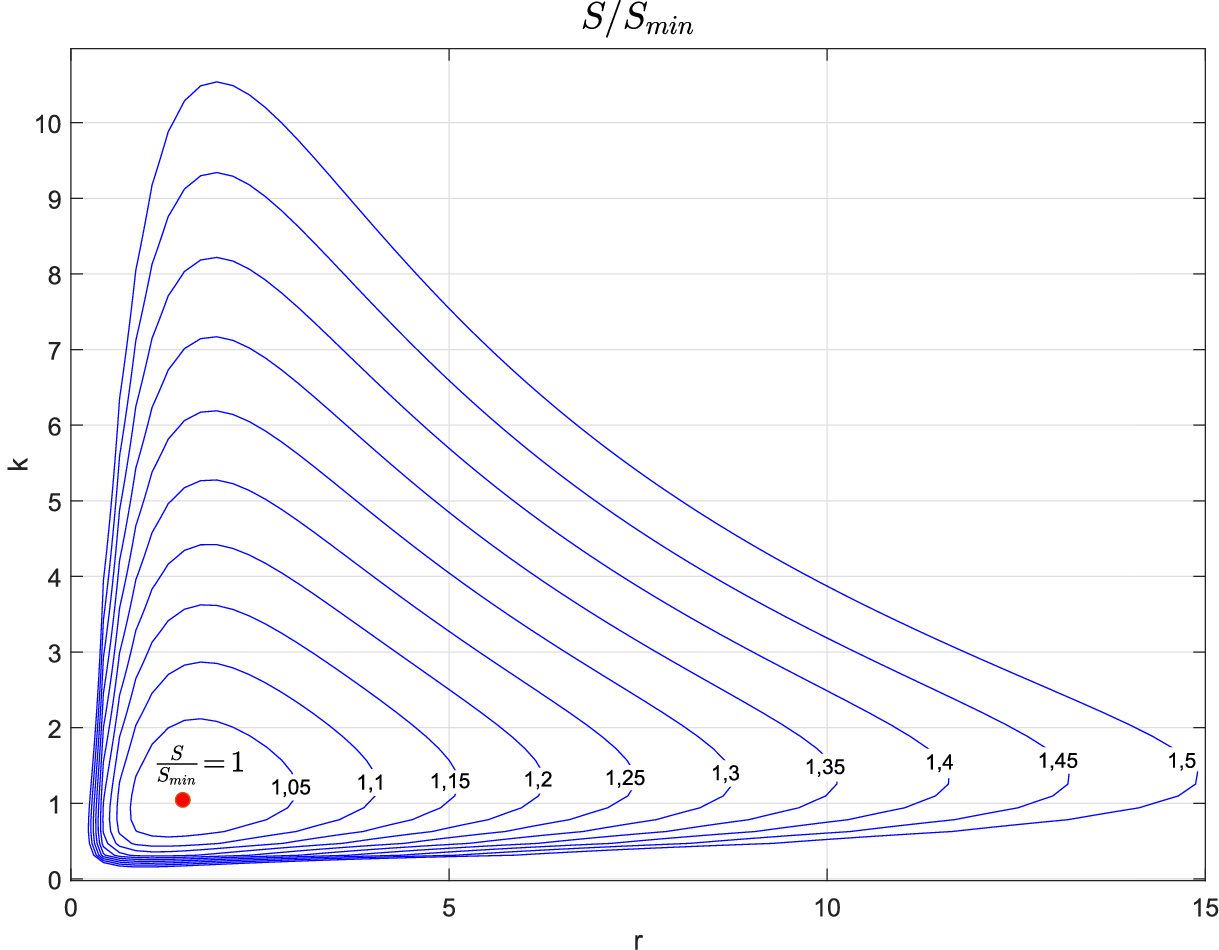}
\caption{Contour plot of the compactness measure $\frac{S}{S_{min}}$ as a function of  footprint aspect ratio $r$ and slenderness aspect ratio $k$. Here $\al=\pi/4=45^\circ$. Red dot is the point corresponding to the global minimum equal to 1. }
\label{poziomice}
\end{figure}

Fig. \ref{poziomice} illustrates the measure $\frac{S}{S_{min}}$
in the form of the contour plot with the same parameters as in Fig. \ref{measure_plot}. The numerical labels on the level curves indicate the corresponding values of $\frac{S}{S_{min}}$. Each level curve thus delineates the pairs $(r,k)$ with the same degree of compactness. Such a graphical representation is of direct utility in the design stage of the barn-house, since it allows for the determination of the building parameters at which the desired compactness is achieved. This makes it possible to determine a house shape that fits well with the building plot, while at the same time its form and compactness remain close to optimal.

\section{Optimal shape for the case of fixed floor area}
\label{area}

In what follows, we consider another optimality problem for the barn-type house, relevant from the practical point of view.
Namely, we assume that the floor area of the house
\[
F = WL
\]
is fixed. For a given $F$ we will find the parameters
$W$ and $L$ such that the external surface of the building $S$ is minimal. In order for this optimization problem to be well-posed, it is necessary to assume that the height
$H$ and the roof slope $\al$ are also fixed. Otherwise, smaller $H$ and $\al$ would always lead to a smaller external surface. It should be emphasized that such a situation is frequently encountered in architectural practice. The investor typically specifies the required floor area of the house as well as the height of the rooms.

In such setting the external surface of the house is given by (cf. \eqref{powierzchnia}):
\begin{align}
S=S(W)&=2WH+2LH+\frac{LW}{\cos(\alpha)}+ W^2 \tan(\alpha)/2 \nonumber \\
&= 2WH+2FH/W+\frac{F}{\cos(\alpha)}+ W^2 \tan(\alpha)/2.
\label{powierzchnia3}
\end{align}
To find the minimum of the function $S(W)$ let us calculate the  derivative $\frac{d S(W)}{d W}$ and solve the equation $\frac{d S(W)}{d W}=0$. We have
\[
\frac{d S(W)}{d W}=2H - 2FH/W^2 + W\tan(\al).
\]
For $W>0$ equation $\frac{d S(W)}{d W}=0$ is thus equivalent to
\eqb
\label{cubic_eq}
W^3\tan(\al)+2HW^2-2FH=0.
\eqe
This is cubic equation of the form
$$
aW^3+bW^2+c=0
$$
with
$$
a=\tan(\al),\;\; b=2H,\;\;c=-2FH.
$$
Using Cardano's formula \cite{tanton2005encyclopedia}
we get that for $H<\tan(\al)\sqrt{27F/16}$ (this condition is certainly satisfied in practice) there exist only one solution of \eqref{cubic_eq}, which determines minimum of external surface $S(W)$. This solution has the form
\begin{align}
W_{min} &=
    -\frac{b}{3a}
    + \frac{2^{1/3}  b^2}{ 3a  \left(-2b^3 - 27a^2c + 3\sqrt{3}\sqrt{4a^2 b^3 c + 27a^4 c^2}\right)^{1/3} } \nonumber \\
    &+ \frac{ \left(-2b^3 - 27a^2c + 3\sqrt{3}\sqrt{4a^2b^3c + 27a^4c^2}\right)^{1/3} }{3a\cdot2^{1/3} } .
\label{WF_min}
\end{align}

In practice, for specified floor area $F$ of the barn house as well as the height of its rooms $H$and angle $\al$,
the above formula determines the optimal width parameter $W_{min}$ as well as the corresponding optimal length parameter $L_{min}=F/W_{min}$ such that the external surface of the building is minimal.

\begin{figure}[t]
\centering
\includegraphics[width=12cm]{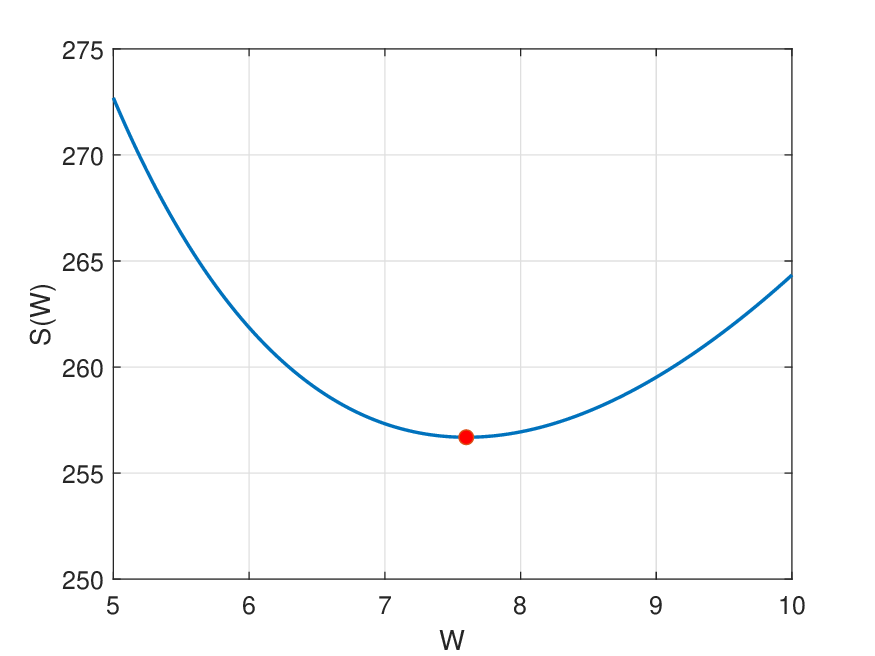}
\caption{Plot of the external surface function $S(W)$  given by \eqref{powierzchnia3} with parameters $F=100 \;m^2$, $H=3\; m$ and $\al=\pi/6=30^\circ$.  The function reaches its minimum at the point $W_{min}=7,60\; m$, the corresponding minimal surface is equal to $S_{min}= 256.69 \;m^2$ (red dot). }
\label{S_W}
\end{figure}

In example, if we take $F=100 \;m^2$, $H=3\; m$ and $\al=\pi/6=30^\circ$, then from \eqref{WF_min} we get that the optimal parameters are $W_{min}=7,60 \;m$ and
$L_{min}=F/W_{min}=13,16 \;m$. In Fig. \ref{S_W} we show the graph of the corresponding external surface function $S(W)$ given by \eqref{powierzchnia3}. Indeed, it reaches its minimum at the point $W_{min}=7,60 \;m$, the corresponding minimal surface is equal to $S_{min}= 256.69 \;m^2$.

Again, we underline that from a practical point of view, the parameters obtained above are perfectly suited for real-life applications in house construction.

\section{Case studies}
\label{cases}
\begin{figure}[t]
\centering
\includegraphics[width=14 cm]{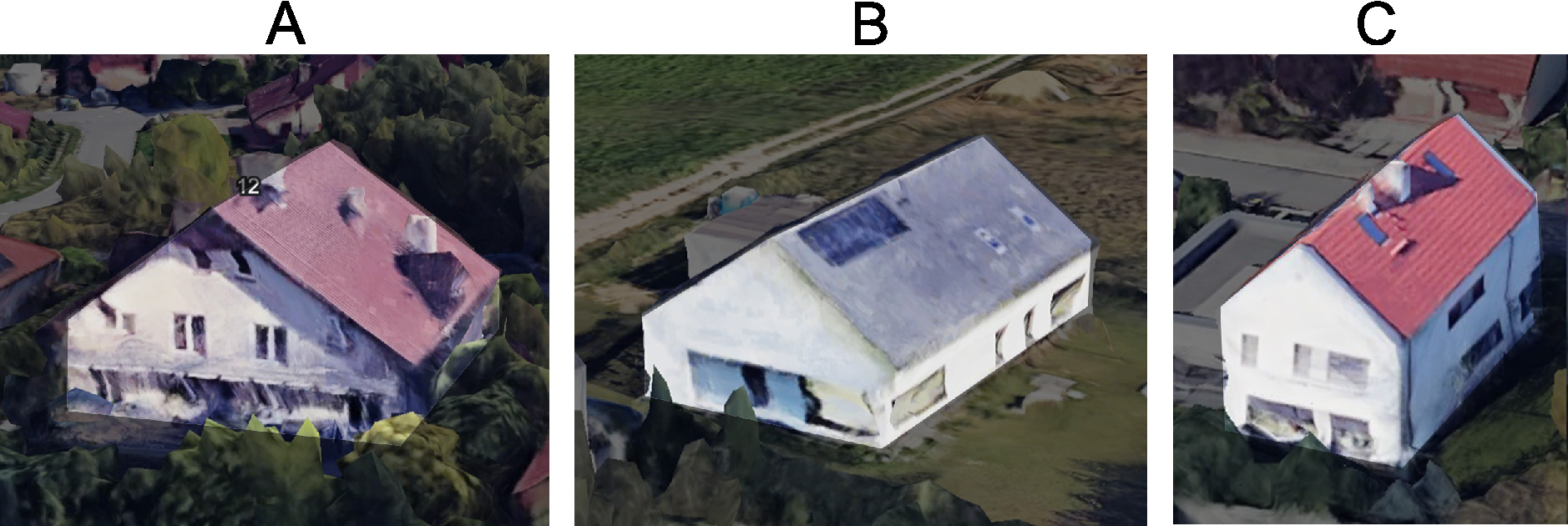}
\caption{Images of the three houses analyzed as case studies using the results obtained in Sections \ref{analiza}--\ref{area}. Images source: Google Earth. }
\label{domy}
\end{figure}

To illustrate the practical applicability of the results developed in previous sections, we analyze three different buildings as case studies. To demonstrate the use of the derived formulas, we analyze three recently built barn-type houses of varying shapes: A, B and C (see Fig. \ref{domy}). The buildings were selected in Wroclaw, Poland, to represent a diverse range of forms and sizes. The images and geometric properties of the houses were derived from Google Earth.

\subsection{Analysis of House A}
The parameters of House A are the following:
width $W=19,9\;m$,
length $L=15,75\;m$,
height $H=5\;m$,
roof slope angle $\al=35^\circ$.
This gives the volume $V = 1567,1\; m^3$ and the external surface
$S = 877,76 m^2$.

Let us first apply the optimality analysis derived in Sec. \ref{analiza}.
For given fixed $V = 1567,1\; m^3$ and $\al=35^\circ$, from Eqs. \eqref{S_min}--\eqref{H_min} we get that the optimal parameters of the house
are: $S_{min} = 737,58 \; m^2$.
$W_{min} = 10,9 \;m$,
$L_{min} = 15,34 \; m$,
$H_{min} = 9,36 \; m$

Consequently, the compactness measure equals
\[
\frac{S}{S_{min}}= 1,19.
\]
The above optimization procedure could reduce the area of the external surface by $S-S_{min}=140,18 \; m^2$. This would result in significant savings in construction costs and energy consumption of the house.

Next, we apply the optimality analysis derived in Sec. \ref{area}.
For given fixed floor area $F = WL = 313,42\; m^2$, height $H=5\;m$ and angle $\al=35^\circ$,
from Eq. \eqref{WF_min} and using the relationship $L_{min}=F/W_{min}$, we get that the optimal parameters of the house
are:
$W_{min} = 12,84 \;m$,
$L_{min} = 24,39 \; m$.
It gives the optimal external surface  $S_{min} = 812,84 \; m^2$ and the ratio $S/S_{min}=1,08 $

The above optimization procedure could reduce the area of the external surface by $S-S_{min}= 64,92 \; m^2$. This would result in considerable savings in construction costs and energy consumption of the house.

\subsection{Analysis of House B}
The parameters of House B are as follows:
width $W=12,4\;m$,
length $L=20,5\;m$,
height $H=4,1\;m$,
roof slope angle $\al=35^\circ$.
This gives the volume $V = 1042,2\; m^3$ and the external surface
$S = 633,93 m^2$.

Let us first apply the optimality analysis derived in Sec. \ref{analiza}.
For given fixed $V = 1042,2\; m^3$ and $\al=35^\circ$, from Eqs. \eqref{S_min}--\eqref{H_min} we get that the optimal parameters of the house
are: $S_{min} = 561,97 \; m^2$.
$W_{min} = 9,52 \;m$,
$L_{min} = 13,39 \; m$,
$H_{min} = 8,18 \; m$

Therefore, the compactness measure equals
\[
\frac{S}{S_{min}}= 1,13.
\]
The above optimization procedure could reduce the area of the external surface by $S-S_{min}=71,96 \; m^2$. This would result in considerable savings in construction costs and energy consumption of the house.

Next, let us apply the optimality analysis derived in Sec. \ref{area}.
For given fixed floor area $F = WL = 254,2\; m^2$, height $H=4,1\;m$ and angle $\al=35^\circ$,
from Eq. \eqref{WF_min} and using the relationship $L_{min}=F/W_{min}$, we get that the optimal parameters of the house
are:
$W_{min} = 11,36 \;m$,
$L_{min} = 22,37 \; m$.
It gives the optimal external surface  $S_{min} = 632,14 \; m^2$ and the ratio $S/S_{min}=1,0028 $

The above optimization procedure could reduce the area of the external surface by $S-S_{min}= 1,79 \; m^2$. This type of analysis shows that House B is very close to optimal.

\subsection{Analysis of House C}
The parameters of House C are the following:
width $W=8\;m$,
length $L=13,5\;m$,
height $H=5,8\;m$,
roof slope angle $\al=40^\circ$.
This gives the volume $V = 626,4\; m^3$ and the external surface
$S = 417,23 m^2$.

Let us first apply the optimality analysis derived in Sec. \ref{analiza}.
For given fixed $V = 626,4\; m^3$ and $\al=40^\circ$, from Eqs. \eqref{S_min}--\eqref{H_min} we get that the optimal parameters of the house
are: $S_{min} = 412,01 \; m^2$.
$W_{min} = 7,72 \;m$,
$L_{min} = 11,15 \; m$,
$H_{min} = 7,28 \; m$

Therefore, the compactness measure equals
\[
\frac{S}{S_{min}}= 1,01.
\]
The above optimization procedure could reduce the area of the external surface by $S-S_{min}=5,22 \; m^2$. This type of analysis shows that House C is very close to optimal.

Now, let us apply the optimality analysis derived in Sec. \ref{area}.
For given fixed floor area $F = WL = 108\; m^2$, height $H=5,8\;m$ and angle $\al=40^\circ$,
from Eq. \eqref{WF_min} and using the relationship $L_{min}=F/W_{min}$, we get that the optimal parameters of the house
are:
$W_{min} = 8,23 \;m$,
$L_{min} = 13,12 \; m$.
It gives the optimal external surface  $S_{min} = 417,09 \; m^2$ and the ratio $S/S_{min}=1,0003 $

The above optimization procedure could reduce the area of the external surface by $S-S_{min}= 0,14 \; m^2$. Again, the analysis shows that House C almost optimal.

Summary of the optimization procedure from Sec. \ref{analiza} (fixed volume) for the three case studied barn-type houses can be found in Tab. \ref{tab:volume}. Summary for the case of fixed floor area and height is presented in Tab. \ref{tab:area}.

\begin{table}[h!]
\scriptsize
\centering
\caption{Summary of the optimization procedure from Sec. \ref{analiza} (fixed volume) for the three case studied barn-type houses.}
\label{tab:volume}
\begin{tabular}{|l|c|c|c|c|c|c|c|c|c|c|c|c|}
\hline
& \multicolumn{6}{c|}{\textbf{Real parameters}} & \multicolumn{4}{c|}{\textbf{Optimal parameters}} & \multicolumn{2}{c|}{\textbf{Comparison}} \\
\hline
& \makecell{$W$\\$[m]$} & \makecell{$L$\\$[m]$} & \makecell{$H$\\$[m]$ } & \makecell{$\al$\\$[^\circ]$ } & \makecell{$V$\\ $[m^3]$} & \makecell{$S$\\ $[m^2]$} & \makecell{$W_{min}$\\ $[m]$} & \makecell{$L_{min}$\\$[m]$} & \makecell{$H_{min}$\\ $[m]$} & \makecell{$S_{min}$\\ $[m^2]$} & \makecell{$S/S_{min}$} & \makecell{$S-S_{min}$\\ $[m^2]$} \\
\hline
\textbf{House A} & 19,9 & 15,75 & 5& 35 & 1567,1 & 877,76 & 10,9 & 15,34 & 9,36 & 737,58 & 1,19 & 140,18 \\
\hline
\textbf{House B} &12,4 &20,5 &4,1 &35 &1042,2 &633,93 &9,52 &13,39 &8,18 &561,97 &1,13 &71,96 \\
\hline
\textbf{House C} &8 &13,5 &5,8 &40 &626,4 &417,23 &7,72 &11,15 &7,28 &412,01 &1,01 &5,22 \\
\hline
\end{tabular}
\end{table}

\begin{table}[h!]
\scriptsize
\centering
\caption{Summary of the optimization procedure from Sec. \ref{area} (fixed floor area and height) for the three case studied barn-type houses.}
\label{tab:area}
\begin{tabular}{|l|c|c|c|c|c|c|c|c|c|c|c|}
\hline
& \multicolumn{6}{c|}{\textbf{Real parameters}} & \multicolumn{3}{c|}{\textbf{Optimal parameters}} & \multicolumn{2}{c|}{\textbf{Comparison}} \\
\hline
& \makecell{$W$\\$[m]$} & \makecell{$L$\\$[m]$} & \makecell{$H$\\$[m]$ } & \makecell{$\al$\\$[^\circ]$ } & \makecell{$F$\\ $[m^2]$} & \makecell{$S$\\ $[m^2]$} & \makecell{$W_{min}$\\ $[m]$} & \makecell{$L_{min}$\\$[m]$} & \makecell{$S_{min}$\\ $[m^2]$} & \makecell{$S/S_{min}$} & \makecell{$S-S_{min}$\\ $[m^2]$} \\
\hline
\textbf{House A} & 19,9 & 15,75 & 5& 35 & 313,42 & 877,76 & 12,84 & 24,39 & 812,84 & 1,08 & 64,92 \\
\hline
\textbf{House B} & 12,4 &20,5 &4,1 &35 & 254,2&633,93 &11,36 &22,37 &632,14 & 1,0028&1,79 \\
\hline
\textbf{House C} &8  &13,5 &5,8 &40 &108 &417,23 & 8,23&13,12 &417,09 &1,0003 &0,14 \\
\hline
\end{tabular}
\end{table}

\subsection{Discussion}
\label{discussion}
The case studies demonstrate how the theoretical framework developed in Sections 2â€“4 translates into practical applications for real buildings. For each analyzed barn-type house, the comparison between real and optimal parameters highlights the degree of compactness and the potential for improvement in terms of reducing the external surface area.

House A shows the largest deviation from the theoretical optimum, with a compactness ratio of 1.19 (fixed volume) and 1.08 (fixed floor area). This indicates that significant reductions in envelope areaâ€”up to 140 mÂ˛â€”could be achieved through adjustments to proportions. Such improvements would directly lower construction costs and reduce heat loss through the envelope, thus enhancing energy efficiency.

House B presents an intermediate case. While its compactness ratio is 1.13 for fixed volume, under the fixed floor area condition it nearly coincides with the theoretical optimum (S/Smin = 1.0028). This suggests that although its proportions could be refined for slightly better performance, the current form already achieves a very good balance between usable space, aesthetics, and energy efficiency.

House C is almost perfectly aligned with the optimal geometry, with compactness ratios of 1.01 and 1.0003 under the two scenarios. The minimal deviations indicate that this building already embodies the design principles suggested by the optimization models. In practice, this means that for certain forms, especially more slender proportions at higher roof slopes, contemporary architectural practice naturally converges toward the theoretical optimum.

Overall, the case studies underline the practical relevance of the  results obtained in Sections \ref{analiza}--\ref{area}. They confirm that theoretical optimization can point to meaningful improvements in design, but they also reveal that many real-world barn-type houses are already close to optimal. This convergence suggests that architectural intuition, plot constraints, and aesthetic preferences often lead to solutions that balance functionality and compactness effectively. The findings also emphasize the potential of the compactness measure as a decision-making tool at the early design stage, allowing architects to quantify trade-offs between form, surface area, and energy use.

\subsection{Software}
\label{software}

In order to enable and facilitate the practical application of the findings presented in this study, we have developed two free and openly accessible computer programs designed to perform optimization analyses of barn-type houses.

\begin{figure}[t]
\centering
\includegraphics[width=16 cm]{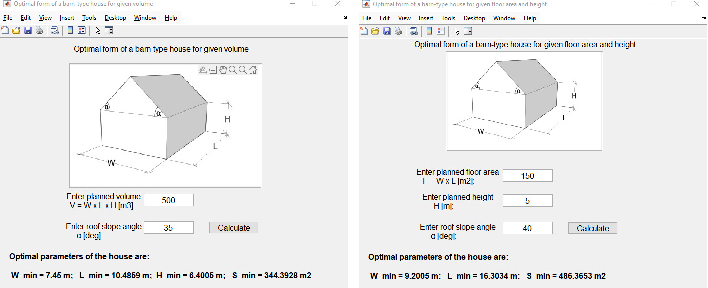}
\caption{Screenshots of two developed and shared computer programs designed to perform optimization analyses of barn-type houses. }
\label{programy}
\end{figure}

The first program performs the analysis of Sec. \ref{analiza} for the case of fixed volume and roof slope angle. The program is available at the web site of Rokita-Projekt architectural Office, at the URL
\url{https://rokita-projekt.pl/barn_house_optimization_fixed_volume.zip} .
Screenshot from the running program is shown if Fig. \ref{programy} (left panel).

The second program performs the analysis of Sec. \ref{area} for the case of fixed floor area, height and roof slope angle. The program is available at the URL \url{https://rokita-projekt.pl/barn_house_optimization_fixed_floor_area.zip} .
Screenshot from the running program is shown if Fig. \ref{programy} (right panel).

Both programs were written in Matlab R2024b environment. MATLAB Runtime (R2024b) package is required to run the programs. The package installer is freely available on the official MathWorks website.

\section{Conclusions}
\label{conclusions}

This study has examined the problem of determining the optimal geometric
parameters of barn-type houses with the aim of minimizing the external
surface of the building envelope while maintaining either a prescribed volume
or a fixed floor area. By formulating the problem mathematically, explicit
formulas for the optimal dimensions of width, length, and height were
derived for given roof slope angles, together with a compactness measure
that provides a dimensionless and scale-independent assessment of how
close a given form is to the theoretical optimum. These results, validated
through both analytical derivations and practical case studies, offer a robust
framework for improving the design of barn-type houses in ways that
enhance energy performance and reduce construction costs.

The findings show that for fixed volume conditions there exist explicit
relationships between footprint aspect ratio and slenderness aspect ratio that
determine the optimal form, leading directly to closed-form expressions for
the minimal surface area of the building. The compactness measure,
expressed as the ratio of the actual surface to the minimal surface, was shown
to be highly effective in quantifying the deviation of real designs from the
optimum. Under conditions of fixed floor area and prescribed height, the
optimization problem reduces to solving a cubic equation, which yields the
optimal width and corresponding length. This situation reflects common
architectural practice, where investors often specify the floor area and room
height, and therefore the results are directly transferable to real-world
applications. The case studies conducted for three barn-type houses in
WrocĹ‚aw demonstrate that the theoretical framework not only provides
valuable insights but also leads to tangible improvements. Some existing
houses deviate significantly from the optimum, and their envelope area could
be reduced considerably, resulting in lower construction costs and improved
thermal efficiency. Other houses were shown to be very close to the
theoretical optimum, which suggests that contemporary architectural
intuition and practice already align closely with the principles identified by
the mathematical optimization.

The implications of the study extend beyond the immediate problem of
geometric compactness. While minimizing the external envelope area
provides clear benefits in terms of material efficiency and energy demand,
future research should aim to integrate the presented methodology with
dynamic energy simulations that take into account climatic conditions,
orientation, and solar exposure. In this way, optimality could be assessed not
only in abstract geometric terms but also in relation to actual heating and
cooling needs across different climate zones. An additional extension of this
work concerns the integration of material considerations, since the choice of
construction systems influences both the embodied carbon and the life-cycle
performance of the building. A multi-objective optimization that balances
geometry, envelope area, and sustainability of materials would provide a
more complete assessment. Furthermore, the economic dimension could be
introduced by incorporating life-cycle cost analysis, allowing the evaluation
of not only construction cost savings but also operational energy use,
maintenance, and long-term value.

The proposed methodology could also be generalized to other house
typologies beyond barn-type forms. Extending the framework to hip roofs,
gable roofs, or courtyard-based configurations would allow a comparative
study of geometric compactness across a broader range of residential
architecture. Likewise, the influence of urban and regulatory constraints,
such as plot dimensions, daylight access, and neighborhood requirements,
should be integrated into the optimization model to enhance its practical
usefulness. Another promising direction lies in expanding the compactness
measure to encompass additional performance indicators, including
daylighting, mechanical ventilation, air conditioning and indoor thermal comfort, thereby
producing a more holistic sustainability metric.

The above open problems are the subject of our future research.

In conclusion, the work presented here provides a rigorous mathematical
foundation for the rational design of barn-type houses. The closed-form
formulas for optimal dimensions, the compactness measure, and the
illustrative case studies together confirm that the approach is both
theoretically sound and practically relevant. At the same time, the results
highlight that many contemporary designs already approximate the
theoretical optimum, suggesting that architectural intuition and mathematical
optimization converge towards similar solutions. With further development the proposed
methodology has the potential to influence architectural practice more
broadly and to contribute significantly to the development of sustainable design guidelines for single-family housing.

\bibliographystyle{plain} 
\bibliography{bibliography}


\section*{Author contributions statement}
E.R-M.: Writing - original draft, Visualization, Validation, Methodology, Investigation, Data curation, Conceptualization. \\ 
B.G.: Writing - review \& editing, Validation, Conceptualization.\\
M.M.: Writing - original draft, Methodology, Investigation, Formal analysis, Software.

\section*{Funding Declaration}
This research received no specific grant from any funding agency in the public, commercial, or not-for-profit sectors.

\section*{Declaration of competing interest}
The authors declare no competing interest.

\section*{Data availability statement}
Data underlying this article will be made available by the corresponding author (Marcin Magdziarz) upon reasonable request.

\end{document}